\documentclass[10pt,reqno]{amsart}

\usepackage{fullpage,amsfonts,amsmath,amssymb,amscd, youngtab}

\input xypic

\theoremstyle{plain}
\newtheorem*{lem}{Lemma}

\newtheorem*{cor}{Corollary}
\newtheorem*{thm}{Theorem}

\newtheorem*{hypo}{Hypothesis}
\theoremstyle{Definition}

\newtheorem*{defn}{Definition}

\newtheorem*{conj}{Conjecture}

\newcommand{\Dim}{\mathrm{Dim}}
\newcommand{\Irr}{\mathrm{Irr}}

\newcommand{\N}{\mathbb{N}}
\newcommand{\Z}{\mathbb{Z}}

\newcommand{\C}{\mathbb{C}}
\newcommand{\Q}{\mathbb{Q}}
\newcommand{\K}{\mathbf{k}}
\newcommand{\blambda}{\boldsymbol{\lambda}}
\newcommand{\bmu}{\boldsymbol{\mu}}
\renewcommand{\leq}{\leqslant}
\renewcommand{\geq}{\geqslant}
\newcommand{\bep}{\boldsymbol{\epsilon}}

\newcommand{\h}{\mathfrak{h}}

\newcommand{\rat}{\overline{H}_{\mathbf{k}}}

\newcommand{\dnorm}{\mathbf{m}}

\newcommand{\one}{\scriptstyle{0}}
\newcommand{\x}{\scriptstyle{1}}
\newcommand{\xx}{\scriptstyle{2}}
\newcommand{\xxx}{\scriptstyle{3}}
\newcommand{\xxxx}{\scriptstyle{4}}
\newcommand{\y}{\scriptstyle{{-1}}}
\newcommand{\yy}{\scriptstyle{{-2}}}

\newcommand{\Res}{\mathrm{Res}}

\title{The Calogero-Moser partition and Rouquier families for complex reflection groups}
\author{Maurizio Martino}

\address{Mathematisches Institut, Universit\"at Bonn, Beringstr. 1, }

\email{mmartino@math.uni-bonn.de}

\begin{document}

\begin{abstract} Let $W$ be a complex reflection group. We formulate a conjecture
relating blocks of the corresponding restricted rational Cherednik algebras and Rouquier families for cyclotomic Hecke algebras.
We verify the conjecture in the case that $W$ is a wreath product of a symmetric group with a cyclic
group of order $l$.
\end{abstract}

\maketitle

\section{Introduction}

\subsection{} Let $W$ be a complex reflection group. The aim of this note is to state,
and for an infinite family of complex reflection groups, prove a
conjecture relating restricted rational Cherednik algebras and
cyclotomic Hecke algebras for complex reflection groups. The former
algebra is finite dimensional factor algebra of a rational Cherednik
algebra with interesting properties, which has been used, for
example, to study the existence of symplectic resolutions of quotient singularities.
Its simple modules are labeled naturally by the set, $\Irr W$, of
simple $\C W$-modules. We can partition $\Irr W$ according to the
blocks of the restricted rational Cherednik algebra; we call this
partition the Calogero-Moser partition (the spectra of the centres
of rational Cherednik algebras are called Calogero-Moser spaces).

\subsection{} Cyclotomic Hecke algebras for complex reflection groups
are objects are current interest, which are expected to provide
insight into the representation theory of finite reductive groups
and to display behaviour analogous to Hecke algebras associated to
series of reductive algebraic groups. In this latter direction,
Rouquier has defined a partition of the set $\Irr W$, which
generalises the notion of Lusztig's families for Weyl groups,
\cite{R1}. This is the partition of $\Irr W$ into Rouquier families.
Conjecture \ref{families} states that the partition into Rouquier
families refines the Calogero-Moser partition and proposes further
numerical connections between them. We expect that this conjecture
is a natural extension of \cite[Conjecture 1.3]{GoMar} which relates
the Calogero-Moser partition to a partition arising (conjecturally)
from cells at unequal parameters for Weyl groups. In the classical
situation this latter partition equals the partition into families,
but since there exists (at present) no cell theory for complex
reflection groups which are not Weyl groups, our conjecture seems to
be the most suitable generalisation of \cite{GoMar}. We prove
Conjecture \ref{families}(i) when $W$ is the wreath product
$G(l,1,n)=\Z/l\Z\wr S_n$ by comparing known combinatorial
descriptions of the two partitions. It would be very interesting to
have a more conceptual understanding of this result.

\subsection{} In section 2 we introduce the main protagonists and state the precise conjecture, which includes
a geometric interpretation of the size of each Rouquier block. In section 3 we prove the conjecture for wreath products
and conclude with an interpretation of the combinatorics via higher level Fock spaces.

\medskip

\noindent \textit{Acknowledgements} The author thanks Iain Gordon for useful comments and
 whose collaboration with the author in \cite{GoMar} was the inspiration for the present work. The author also thanks Gwyn
Bellamy, Maria Chlouveraki and Nicolas Jacon for useful comments.
This work was supported by the SFB/TR 45 ``Periods, Moduli Spaces
and Arithmetic of Algebraic Varieties" of the DFG (German Research
Foundation).

\section{Calogero-Moser partition, Rouquier families and the conjecture}

\subsection{Notation}\label{notation} Let $W$ be a complex reflection group and $\h$ its
reflection representation over $\C$. %Let $\mathrm{Ref}W$ denote the
%set of reflections in $W$ and
Let $\mathcal{A}$ be the set of reflecting hyperplanes in $\h$. Let
$\h^{\mathrm{reg}} = \h \setminus (\bigcup_{H \in \mathcal{A}} H)$.
Given a hyperplane $H \in \mathcal{A}$ we define $W_H$ to be the
subgroup of $W$ of elements that fix $H$ pointwise. Let $e_H =
|W_H|$, and for $\mathcal{C} \in \mathcal{A}/W$ let
$e_{\mathcal{C}}$ be the common value $e_H$ for $H\in \mathcal{C}$.
For every $H \in \mathcal{A}$ we choose $v_H \in \h$ such that $\C
v_H$ is a $W_H$-stable complement to $H$, and choose also a linear
form $\alpha_H \in \h^*$ with kernel $H$. Let $<\ ,\ >$ denote the
natural pairing of $\h^*$ with $\h$. Let $A$ be a finite dimensional algebra, 
then we denote by $\Irr A$ the set of irreducible $A$-modules. We will write $\Irr W$ for
 the set $\Irr \C W$.

%We use the convention $\otimes = \otimes_{\C}$.

\subsection{Rational Cherednik algebras}\label{RCAdefn} We introduce parameters
$\K:=(k_{\mathcal{C},i})_{\mathcal{C} \in \mathcal{A}/W,\ 0\leq i
\leq e_{\mathcal{C}}-1}$ where $k_{\mathcal{C},i} \in \C$ for all
$\mathcal{C}$ and $i$, and $k_{\mathcal{C},0} = 0$ for any
$\mathcal{C}$. We use the convention here and throughout this paper
that the subscript $j$ in $k_{\mathcal{C},j}$ is considered modulo
$e_{\mathcal{C}}$.

The \textit{rational Cherednik algebra} (at $t=0$), $H_{\K}$, is
the quotient of $T(\h \oplus \h^*)\rtimes W$, the smash
product of the free $\C$-algebra on $\h \oplus \h^*$ with $W$, by
the relations:
\begin{align*} [x,x']=0,\ [y,y']=0\ \mathrm{and}\ [y,x] = \sum_{H \in
\mathcal{A}} \frac{<\alpha_H,y><x,v_H>}{<\alpha_H,v_H>} \gamma_H
\end{align*}
for $y,y'\in \h$ and $x,x'\in \h^*$. Here we set \[\gamma_H = \sum_{w\in
W_H\setminus \{1\}} \left( \sum_{j=0}^{e_H - 1} \mathrm{det}(w)^{j}\\
(k_{\mathcal{C},j+1} - k_{\mathcal{C},j}) \right) w\] for all $H\in
\mathcal{A}.$

This is the definition given in \cite{GGOR} and, as noted there, is
equivalent to that given in \cite{EG}. It follows
easily from the definition that \begin{equation}\label{rescale}
H_{\K} \cong H_{ \lambda \K}\end{equation} for any
$\lambda \in \C^*$.

\subsection{The Calogero-Moser partition}\label{CMpart} By \cite[Proposition 4.15]{EG}, there is an algebra
embedding of $A:= \C[\h]^W \otimes_{\C} \C[\h^*]^W$ into the centre of
$H_{\K}$. Let $A_+ \subset A$ denote the maximal ideal of
polynomials with positive degree. The \textit{restricted rational
Cherednik algebra}, $\overline{H}_{\K}$, is the finite dimensional
factor algebra $H_{\K}/A_+H_{\K}$. For more details on the structure
of $\overline{H}_{\K}$, see \cite{Go1}.

The \textit{Calogero-Moser partition}, or $CM_{\K}$-partition, is
defined by the equivalence relation on $\Irr \rat$ given by: $M
\sim_{CM_{\K}} N$ if and only if $M$ and $N$ lie in the same block
of $\rat$. The set of irreducible $\rat$-modules, $\Irr \rat$, can
be identified with the set $\Irr W$, \cite[Proposition 4.3]{Go1},
and so we think of the $CM_{\K}$-partition as a partition of $\Irr
W$. The proof of the next lemma follows directly from
(\ref{rescale})
\begin{lem}
Let $\K$ be a parameter as in \ref{RCAdefn} and let $\lambda \in
\C$. Then the $CM_{\K}$-partition and $CM_{\lambda \K}$-partition
are equal.
\end{lem}

\subsection{Geometric interpretation}\label{CMfibres} Let
$Z_{\K}$ be the centre of $H_{\K}$. The embedding
\[ A \hookrightarrow Z_{\K}\] induces a morphism of schemes
\[ \Psi: \mathrm{Spec}\ Z_{\K} \to \h/W \times \h^*/W.\] We write
$\Psi^*(0)$ for the scheme theoretic fibre of $0$ and $\Psi^{-1}(0)$
for the closed points in $\Psi^*(0)$. By \cite[Corollary 5.8]{Go1}
the $CM_{\K}$-partition is trivial (that is, each equivalence class
is a singleton set) if and only if $\Psi^{-1}(0)$ consists of smooth
points, and this occurs if and only if all irreducible
$\overline{H}_{\K}$-modules have dimension $|W|$.

\subsection{Generic Hecke algebras}\label{genericHecke}
We retain the notation of \ref{notation}. For every $d
> 1$, set $\eta_d = e^{\frac{2\pi \sqrt{-1}}{d}}$ and denote by
$\mu_d$ the group of all $d$th roots of unity. Let $\mu_{\infty}$ be
the group of all roots of unity in $\C$ and let $K$ be a number
field contained in $\Q (\mu_{\infty})$ such that $K$ contains
$\mu_{e_{\mathcal{C}}}$, for all $\mathcal{C}\in \mathcal{A}/W$. We
denote by $\mu(K)$ the group of roots of unity in $K$ and by $\Z_K$
the ring of integers in $K$.

Let $x_0 \in \h^{\mathrm{reg}}$ and denote by $\overline{x}_0$ its
image in $\h^{\mathrm{reg}}/W$. Let $B$ be the fundamental group
$\Pi_1(\h^{\mathrm{reg}}/W,\overline{x}_0)$. Let $\mathbf{u} =
(u_{\mathcal{C},j})_{\mathcal{C}\in \mathcal{A}/W,\ 0\leq j \leq
e_{\mathcal{C}}-1}$ be a set of indeterminates, and let $\Z
[\mathbf{u} , \mathbf{u}^{-1}] := \Z [u_{\mathcal{C},j}^{\pm 1}:
\mathcal{C}\in \mathcal{A}/W,\ 0\leq j \leq e_{\mathcal{C}}-1]$. The
\textit{generic Hecke algebra}, $\mathcal{H}_W$, is the quotient of
$\Z [\mathbf{u} , \mathbf{u}^{-1}] B$ by relations of the form \[
(\mathbf{s}-u_{\mathcal{C},0})(\mathbf{s}-u_{\mathcal{C},1}) \dots
(\mathbf{s}-u_{\mathcal{C},e_{\mathcal{C}-1}}), \] where
$\mathcal{C} \in \mathcal{A}/W$ and $\mathbf{s}$ runs over the set
of monodromy generators around the images in $\h^{\mathrm{reg}}/W$
of the hyperplane orbit $\mathcal{C}$, see \cite[\S 4.C]{BMR}. We
shall from now on assume the following.

\begin{hypo}
$\mathcal{H}_W$ is a free $\Z [\bf{u} , \bf{u}^{-1}]$-module of rank
$|W|$ and $\mathcal{H}_W$ has a symmetrising form $t: \mathcal{H}_W
\to \Z [\bf{u} , \bf{u}^{-1}]$ that reduces to the standard
symmetrising form on $\Z_K W$ upon specialising $u_{\mathcal{C},j}$ to
$\eta_{e_{\mathcal{C}}}^j$. Furthermore, let $\mathbf{v} =
(v_{\mathcal{C},j})_{\mathcal{C}\in \mathcal{A}/W,\ 0\leq j \leq
e_{\mathcal{C}}-1}$ be indeterminates such that $u_{\mathcal{C},j}=
\eta_{e_{\mathcal{C}}}^j v_{\mathcal{C},j}^{|\mu(K)|}$ for all
$\mathcal{C}, j$: then the $K(\mathbf{v})$-algebra $K(\mathbf{v})
\mathcal{H}_W$ is split semisimple.
\end{hypo}

%We have cheated slightly here: normally ones adds an extra condition
%on the form $t$, \cite[2.1]{BK}, and then the final condition
%follows
 %from the other ones, \cite[5.2]{Ma}.
It is known that all but a finite number of complex reflection
groups satisfy this hypothesis, see \cite{AK}, \cite{BM} and
\cite{Ar2}, and it is conjectured to hold for all complex reflection
groups. By Tits' deformation theorem (see, for example,
\cite[Theorem 7.2]{GP}) the specialisation $v_{\mathcal{C},j}
\mapsto 1$ induces a bijection $\Irr W \to \Irr K(\mathbf{v})
\mathcal{H}_W;\ S \mapsto S^{K(\mathbf{v})}$ such that $\Dim_{\C} S
= \Dim_{K(\mathbf{v})} S^{K(\mathbf{v})}$.

\subsection{Cyclotomic Hecke algebras}\label{cyclotomic}

\begin{defn}
A cyclotomic Hecke algebra is the $\Z_K [y,y^{-1}]$-algebra induced
from $\Z_K [\mathbf{v},\mathbf{v}^{-1}] \mathcal{H}_W$ by an algebra
homomorphism of the form \begin{align*} \Z_K
[\mathbf{v},\mathbf{v}^{-1}]& \to \Z_K [y,y^{-1}]\\
v_{\mathcal{C},j} &\mapsto y^{m_{\mathcal{C},j}} \end{align*} such
that \begin{itemize}\item[(i)] $m_{\mathcal{C},j}\in \Z$ for all
$\mathcal{C}\in \mathcal{A}/W$ and $0 \leq j \leq
e_{\mathcal{C}}-1$; \item[(ii)] Set $x:= y^{|\mu (K)|}$. If $z$ is
an indeterminate then the element of $\Z_K [y,z]$ defined by
\[ \Gamma_{\mathcal{C}}(y,z) = \prod_{j=0}^{e_{\mathcal{C}} -1} (z -
\eta_{e_{\mathcal{C}}}^j y^{m_{\mathcal{C},j}})\] is invariant under
$\mathrm{Gal}(K(y)/K(x))$ for all $\mathcal{C}\in \mathcal{A}/W$. In
other words, $\Gamma_{\mathcal{C}}(y,z)$ is contained in $\Z_K
[x^{\pm 1}, z]$.
\end{itemize}
We write $\mathbf{m} = (m_{\mathcal{C},j})_{\mathcal{C}\in
\mathcal{A}/W,\ 0 \leq j \leq e_{\mathcal{C}}-1}$ and denote this
algebra by $\mathcal{H}_W({\mathbf{m}})$.
\end{defn}

The algebra $K(y) \mathcal{H}_W(\mathbf{m})$ has a symmetric form
induced by $t$ and is split semisimple by \cite[\S 4.3]{Ch}. Thus by
Tits' deformation theorem we have bijections
\begin{align*} \Irr W \cong \Irr K(y)
\mathcal{H}_W(\mathbf{m}) \cong \Irr K(\mathbf{v})
\mathcal{H}_W.\end{align*}

\subsection{Rouquier families}\label{families}
We define the \textit{Rouquier ring} to be $\mathcal{R}(y) = \Z_K
[y, y^{-1}, (y^n - 1)^{-1}_{n \geq 1}]$. By Hypothesis
\ref{genericHecke}, $\mathcal{R}(y)\mathcal{H}_W(\mathbf{m}) \subset
K(y) \mathcal{H}_W(\mathbf{m})$ is free of rank $|W|$. We define an
equivalence relation on $\Irr K(y) \mathcal{H}_W(\mathbf{m}) = \Irr
W$ by: $S \sim_{R_{\mathbf{m}}} T$ if and only if $S$ and $T$ belong
to the same block of $\mathcal{R}(y)\mathcal{H}_W(\mathbf{m})$. We
call the equivalence classes of this relation \textit{Rouquier
families}, and we call the subalgebras $K(y) B \subseteq K(y) \mathcal{H}_W(\mathbf{m})$, where $B$ is a block of
$\mathcal{R}(y)\mathcal{H}_W(\mathbf{m})$, \textit{Rouquier blocks}.

\begin{conj}
Let $W$ be a complex reflection group satisfying Hypothesis
\ref{genericHecke}. \\ Let $\K=(k_{\mathcal{C},j})_{\mathcal{C} \in
\mathcal{A}/W,\ 0\leq j \leq e_{\mathcal{C}}-1}$ be a parameter as
in \ref{RCAdefn} such that $k_{\mathcal{C},j} \in \Z$ for all
$\mathcal{C}$ and $j$. Let $\mathbf{m} =
(m_{\mathcal{C},j})_{\mathcal{C}\in \mathcal{A}/W,\ 0 \leq j \leq
e_{\mathcal{C}}-1}$ where $m_{\mathcal{C},j} =
-k_{\mathcal{C},e_{\mathcal{C}}-j}$ for all $\mathcal{C}$ and $j$.
Then:
\begin{itemize}\item[(i)] The partition of $\Irr W$
into Rouquier families associated to $\mathcal{H}_W(\mathbf{m})$
refines the $CM_{\K}$-partition. For generic values of $\K$ the partitions are equal;
\item[(ii)] Let $\K$ be a parameter such that the $CM_{\K}$-partition and partition into Rouquier families are equal.
Let $p \in \Psi^{-1}(0)$ and let $K(y) B$ be the corresponding Rouquier block. Then
$\Dim_{\C}\big( \C[\Psi^*(0)_p] \big) = \Dim_{K(y)} K(y)B.$
\end{itemize}
\end{conj}

A priori there is no reason why the specialisation
$m_{\mathcal{C},j} = -k_{\mathcal{C},e_{\mathcal{C}}-j}$ satisfies
condition (ii) of Definition \ref{cyclotomic}. By Lemma \ref{CMpart}
however, we can assume without loss of generality that $|\mu(K)|$
divides all of the $k_{\mathcal{C},j}$. It seems necessary to
include the condition that $\K$ has integer entries since the Hecke
algebra is only defined at integer parameters. If $\K$ has rational
entries then $\lambda \K$ has integer entries for some $\lambda \in
\Z$, so we can state a version of the conjecture for $\K$ by using
Lemma \ref{CMpart}.

\subsection{Evidence for the conjecture}By the Shepherd-Todd classification of complex
reflection groups, \cite{ST}, $W$ is either a member of the infinite
family $G(l,p,n)$ where $l,n \in \N$ and $p$ divides $l$, or one of
$34$
exceptional groups $G_4, \dots, G_{37}$. %Combining the works
%\cite{Bel} and \cite{Ch} one can show when $W \neq G(m,1,n)$ or $G_4$ then both partitions are
%nontrivial for all values of $\K$.

\begin{itemize}
\item[$\bullet$] Combining \cite[Proposition 16.4(ii)]{EG}, \cite[Proposition 7.3]{Go1}
and \cite[Proposition 3.2]{Bel} we have that when $W$ is not
$G(l,1,n)$ or $G_4$ there exist irreducible $\rat$-modules of
dimension $<|W|$ for all values of $\K$. By \ref{CMfibres} this
implies that in these cases the $CM_{\K}$-partition is nontrivial.
The work of \cite{Ch} shows that the Rouquier families of
$\mathcal{H}_W(\dnorm)$ are never trivial in these cases.

\item[$\bullet$] Suppose now that $W = G(l,1,n)$: this is the wreath product $\Z/l\Z \wr S_n $.
We prove part (i) of the conjecture for these groups in the next section, Corollary \ref{mainthm}.
%There is combinatorial description of the $CM_{\K}$-partition given
%in \cite{Go2}, see also \cite[Theorem 2.5]{GoMar}. On the other hand
%there is also a combinatorial description of the corresponding
%Rouquier families from \cite{BK}. It is shown in \cite{Mar} that
%these are the same, confirming the conjecture in this case.

\item[$\bullet$] Let $W = G(l,1,n)$. Then, for generic rational values
of $\K$, $X_{\K}$ is smooth so that the $CM_{\K}$-partition is
trivial. In this case if we take $p \in \Psi^{-1}(0)$ then by
\cite[Corollary 5.8]{Go1}, $\Dim_{\C}\big( \C[\Psi^*(0)_p] \big)=
(\Dim S)^2$, where $S \in \Irr W$ is the module corresponding to
$p$. On the other hand we know from Corollary \ref{mainthm} that the
Rouquier partition is trivial and so the Rouquier blocks of $K(y)
\mathcal{H}_W(\mathbf{m})$ are simply the blocks of this algebra.
Now by Tits' deformation theorem we have that the dimension of the
block corresponding to $S$ has $K(y)$-dimension $(\Dim S)^2$.
\end{itemize}

\section{Proof of the conjecture for $G(l,1,n)$}

\subsection{Parameters for $G(l,1,n)$}\label{parameters} We fix throught positive integers $l$ and $n$.
The wreath product $W= \Z/l\Z \wr S_n$ has two orbits of reflecting
hyperplanes, $\mathcal{C}_1$ and $\mathcal{C}_2$, with
$e_{\mathcal{C}_1}=2$ and $e_{\mathcal{C}_2}=l$. We will use the
parameter set $\mathbf{h}=(h; H_1, \dots ,H_{l-1})$ where $h=
k_{\mathcal{C}_1, 1}$ and $H_i = k_{\mathcal{C}_2, l-i+1} -
k_{\mathcal{C}_2, l-i}$ for each $i$. We set also $H_0 = -H_1- \dots
-H_{l-1}$.

We will assume throughout that $\mathbf{h}$ has rational entries. We
also assume that $h=-1$ and we fix a positive integer $d$ such that
$d\mathbf{h}= (-d; dH_1, \dots ,dH_{l-1})$ has integer entries. %By
%Lemma \ref{CMpart} the $CM_{\mathbf{h}}$-partition equals the
%$CM_{d\mathbf{h}}$-partition.
Our parameter $\bf m$, calculated with respect to $d\mathbf{h}$, is given by
$m_{\mathcal{C}_1,0}=0, m_{\mathcal{C}_1,1}=d$ and
$m_{\mathcal{C}_2,j} = d\sum_{i=1}^j H_i$ for all $0\leq j \leq
l-1$.

\subsection{Sequences} Let $r\in \Z$. A strictly decreasing sequence of integers $C=\{C_1,
C_2, \dots \}$ will be said to \textit{stabilise with respect to r}
if there exists an $K \geq 1$ such that $C_i = r+1-i$ for all $i
\geq K$. Given any strictly decreasing set of integers, $C$, we define its power series
\[ \pi(C):= \sum_{i\geq 1} x^{C_i}.\] For any positive integer $K$ we define the truncated
power series $\pi(C)_{\leq K} := \sum_{i: C_i\leq K} x^{C_i}$. Given
any sequence of integers, $C$, and any $i\in \Z$ we define $S^iC$ to
be the sequence $\{C_1+i, C_2+i, \dots \}$.

The next notion will be useful to us.

\begin{defn} Let $\mathbf{r} = (r_0, \dots ,r_{l-1})\in \Z^l$. Let $\mathbf{C} = (C^q)_{q=0}^{l-1}$ be an $l$-tuple
where each $C^q$ is a strictly decreasing sequence of integers which
stabilises with respect to $r_q$. We define $\chi(\mathbf{C})$ to be the set
\begin{align*} \bigcup_{i=0}^{l-1} \{l(C^q_i-1)+i+1: 0\leq q \leq
l-1, i\geq 1\}.\end{align*}\end{defn}

We note that $\chi$ is one-to-one: for $l$-tuples of strictly
decreasing integers, $\mathbf{C}$ and $\mathbf{D}$,
$\chi(\mathbf{C}) = \chi(\mathbf{D})$ implies $\mathbf{C} =
\mathbf{D}$. We will need the following elementary result.

\begin{lem}\label{sequenceslemma}
Let $\mathbf{C}$ be as above and let $r = \sum_{i=0}^{l-1} r_i$.
Then the set $\chi(\mathbf{C})$ can be rearranged into a strictly
decreasing sequence of integers which stabilises with respect to
$r$.
\end{lem}
\begin{proof}
Let $d$ be the smallest integer such each $C^q$ stabilises after $d$
steps. Let $r_{\rm min} = \mathrm{min}\{r_0, \dots, r_{l-1}\}$. For
each $q$ let ${D}^q$ be the sequence obtained by removing the first
$d + (r_q - r_{\rm min})$ terms from $C^q$. Thus $D^q=\{r_{\rm
min}-d,r_{\rm min}-d-1, \dots \}$. Let $\mathbf{D}=(D^0, \dots
,D^{l-1})$. Then $\chi(\mathbf{D}) = \{l(r_{\rm min}-d),
l(r_{\rm{min}}-d)-1, \dots \}$ and $|\chi(\mathbf{C}) \setminus
\chi(\mathbf{D})| = \sum_{q=0}^{l-1} d + (r_q - r_{\rm min}) = r +
l(d-r_{\rm min})$.
\end{proof}

\subsection{Partitions} A \textit{partition} of $n$ is a sequence of natural numbers
$\lambda = (\lambda_1 \geq \lambda_2 \dots \geq \lambda_k)$ such
that $|\lambda|:= \sum_{i=1}^k \lambda_i = n$. The integer $k$ is
called the \textit{length} of $\lambda$, and we will denote this by
$L(\lambda)$. We use the convention that $\lambda_i = 0$ for $i>k$,
and we denote the set of partitions of $n$ by $\mathcal{P}(n)$. The
\textit{Young diagram} of $\lambda$ is $Y(\lambda) := \{(a,b)\in
\Z^2: 1\leq a \leq L(\lambda),\ 1 \leq b \leq \lambda_a\}$. The
elements of the Young diagram are called \textit{nodes} and we
define the content of a node to be $\mathrm{cont}(a,b) = b-a$. For
example, the Young diagram of the partition $(5,3,2)$ with nodes
labeled with their content is:
\[\young(\one \x \xx \xxx \xxxx,\y \one
\x,\yy \y)\]

\subsection{$\beta$-numbers}\label{betanumbers} Let $r \in \Z$. We define the $r$-shifted $\beta$-number of
$\lambda$ to be the sequence of decreasing integers \begin{align*}
\beta^r(\lambda)=\{\lambda_1+r, \lambda_2+r-1, \dots
,\lambda_j+r+1-j, \dots \}.\end{align*} We define $\beta^r(\lambda)_i
= \lambda_i +r +1 -i$ for each $i \geq 1$. The decreasing
sequence $\beta^r(\lambda)$ stabilises with respect to $r$; conversely, any decreasing
sequence of integers which stabilises with respect to some $r\in \Z$
equals the $r$-shifted $\beta$-number of a unique partition. When $r=0$ we shall
simply write $\beta(\lambda)$ instead of $\beta^r(\lambda)$.

%It will sometimes be useful to represent a $\beta$-number (or indeed
%any decreasing set of integers) by its power series, that is, $\pi
%(\beta^r(\lambda))(x) = \pi (\beta^r(\lambda)) := \sum_{i \geq 1}
%x^{\beta^r(\lambda)_i}$.

The \textit{residue} of a partition $\lambda$, $\mathrm{Res}_{\lambda}(x)$,
is the element of $\Z[x^{\pm 1}]$ given by \[\sum_{(a,b) \in
Y(\lambda)} x^{\mathrm{cont}(a,b)}.\] For $r \in \Q$ we define the
$r$-shifted residue of $\lambda$ to be $\mathrm{Res}^r_{\lambda}(x)
:= x^r \mathrm{Res}_{\lambda}(x)$. %It will be useful to extend this
%definition slightly by allowing $r$ to be a rational number, in
%which case $\Res^r_{\lambda}(x)$ is defined as above and lies in the
%ring $\Z<x^a: a\in \Q>$.
We have the following relationship between residues and
$\beta$-numbers:
\begin{align}\label{res-cont} (x-1)\Res^r_{\lambda}(x) = \pi(\beta^r(\lambda))_{\leq K} -
\frac{x^{r+1-K}-x^{r+1}}{1-x}\end{align} for any $K \geq
L(\lambda)$.

\subsection{$J$-hearts}\label{Jhearts} Given $j \in \{0, 1,
\dots , l-1\}$ we say a node $(a,b)\in Y(\lambda)$ is
$j$\textit{-removable} if $\mathrm{cont}(a,b)$ is equal to $j$
modulo $l$ and if $Y(\lambda) \setminus (a,b)$ is the Young diagram
of some partition. Given a subset $J \subseteq \{0,\dots ,l-1\}$ we
define the $J$-\textit{heart} of $\lambda$ to be the partition
obtained by removing as often as possible $j$-removable boxes, where
$j \in J$. Denote this partition $\lambda_J$.

The notion of $j$-removability is related to $\beta(\lambda)$ as
follows. If a node, $(a,b)$, is removable then it lies at the right
hand edge of row $a$ and so $b = \lambda_a$. Now $(a,\lambda_a)$ is
$j$-removable for some $j \in J$ if and only if its content is
congruent to $j$ modulo $l$ and $(a+1,\lambda_a)$ does not lie in
$Y(\lambda)$. This is equivalent to: $\beta(\lambda)_a - 1 \equiv j\
\mathrm{mod}\ l$ and $\beta(\lambda)_{a+1} < \beta(\lambda)_a - 1$.
Furthermore, the $\beta$-number of $\lambda \setminus (a,\lambda_a)$
is $(\beta(\lambda) \setminus \beta(\lambda)_a) \cup
\{\beta(\lambda)_a-1\}$.

\subsection{Residues and Rouquier families}\label{multiparts} An \textit{$l$-multipartition} of $n$ is an
$l$-tuple $(\lambda^{(0)}, \dots ,\lambda^{(l-1)})$ of partitions
such that $\sum_{i=0}^{l-1} |\lambda^{(i)}| =n$. We denote the set
of $l$-multipartitons of $n$ by $\mathcal{P}(l,n)$. There is a
natural bijection between $\mathcal{P}(l,n)$ and $\Irr W$, see
\cite[6.1.1]{R3}, for example. Thus we think of the
$CM_{\K}$-partition and the partition into Rouquier families as
partitions of the set $\mathcal{P}(l,n)$.

Given $\mathbf{r}= (r_0, \dots r_{l-1}) \in \Z^l$ we define the
$\mathbf{r}$-shifted residue of $\blambda$ to be
\[ \mathrm{Res}^{\mathbf{r}}_{\blambda}(x) := \sum_{i=0}^{l-1}
\mathrm{Res}^{r_i}_{\lambda^{(i)}}(x).\]

\begin{thm}\cite[Th\'{e}or\`{e}me 3.13]{BK}, \cite[Theorem
3.11]{Ch2} Let $\blambda, \bmu \in \mathcal{P}(l,n)$ and let
$\mathbf{m}$ be as in \ref{parameters}. Define
$\overline{\mathbf{m}} = (0, dH_1, dH_1 + dH_2, \dots ,dH_1+\dots
+dH_{l-1})$. If $\blambda \sim_{R_{\mathbf{m}}} \bmu$ then
\begin{equation*}\mathrm{Res}^{\overline{\mathbf{m}}}_{\blambda}(x^d) =
\mathrm{Res}^{\overline{\mathbf{m}}}_{\bmu}(x^d).
\end{equation*}
\end{thm} For generic values of ${\mathbf{m}}$ the converse is also true,
since for generic ${\mathbf{m}}$ the partition into Rouquier
families is trivial, \cite[Proposition 3.12]{Ch2}; in fact, if $l$
is a power of a prime number then the converse statement is true for
all ${\mathbf{m}}$.

\subsection{} Let $\blambda \in \mathcal{P}(l,n)$ and $\mathbf{r}\in \Z^l$.
We obtain a partition from $\blambda$ and $\mathbf{r}$ as
follows. Set $r = \sum_{i=0}^{l-1} r_i$. Let
$\boldsymbol{\beta}^{\bf r} = (\beta^{r_0}(\lambda^{(0)}), \dots
,\beta^{r_{l-1}}(\lambda^{(l-1)}) )$. The set
$\chi(\boldsymbol{\beta}^{\bf r})$ can be arranged into a decreasing
sequence of integers stabilising to $r$, Lemma \ref{sequenceslemma}.
Therefore $\chi(\boldsymbol{\beta}^{\bf r})$ equals
$\beta^r(\tau_{\mathbf{r}}(\blambda))$ for some partition
$\tau_{\mathbf{r}}(\blambda)$, and we obtain a map
\begin{align*} \Z^l \times \coprod_n \mathcal{P} (l,n) \to \coprod_n
\mathcal{P}(n); \ \ (\mathbf{r},\blambda) \mapsto
\tau_{\mathbf{r}}(\blambda).\end{align*}

\subsection{Affine symmetric group}\label{affinesymgroup} Throughout this subsection subscripts are considered modulo $l$.
Let $S_l$ denote the symmetric group on $l$ letters. We
identify $S_l$ with permutations of the set $\{0, \dots ,l-1\}$,
which is generated by elements $s_i$ for $1\leq i \leq l-1$, where
$s_i$ is the simple transposition swapping $i-1$ and $i$. There is an action of $S_l$ on $\Z^l$ via:
\begin{align}\label{affinesymaction} s_i\cdot (\theta_0, \dots
,\theta_{l-1}) = (\theta_0 , \dots , \theta_{i-1} + \theta_i,
-\theta_i, \theta_i + \theta_{i+1}, \dots ,\theta_{l-1})\ \ \
\mathrm{for\ all}\ 1\leq i \leq l-1.\end{align} Let $e_0, \dots ,e_{l-1}$ denote the standard basis of the lattice
$\Z^l$. Let $R$ denote the
root lattice of type $\tilde{A}_l$, which is the sublattice of
$\Z^l$ generated by the simple roots $\alpha_i=-e_{i-1} + 2e_i -
e_{i+1}$ for all $0 \leq i \leq l-1$. The action of $S_l$ preserves
$R$ and we define the affine symmetric group, $\tilde{S}_l$, to be
the semidirect product $R \rtimes S_l$.

The equations in (\ref{affinesymaction}) extend naturally to define
an action of ${S}_l$ on $\Q_1^l := \{ (\theta_0, \dots
,\theta_{l-1})\in \Q^l : \theta_0 + \dots + \theta_{l-1}=1\}$. The
group $R$ acts on $\Q_1^l$ by translations. These two actions
combine to give an action of $\tilde{S}_l$ on $\Q_1^l$. Let
\[\mathcal{A} = \{ (\theta_0, \dots ,\theta_{l-1})\in \Q_1^l: 0
\leq \theta_i \leq 1\ \mathrm{for}\ 0\leq i\leq l-1 \}\] By
\cite[Proposition 4.3]{Hum}, for every $\theta \in \Q_1^l$ there is
a (not necessarily unique) $w_{\theta} \in \tilde{S}_l$ such that
$w_{\theta} \cdot \theta \in \mathcal{A}$. Let $w_{\theta} \cdot
\theta = (\varepsilon_0, \dots ,\varepsilon_{l-1})$. We define the
\textit{type}, $J$, of an element $\theta \in \Q_1^l$ to be the set
$\{j\in \{0, \dots, l-1\}: \varepsilon_j=0\}$.

\subsection{$CM_{\mathbf{h}}$-partition}\label{Iainthm} We describe
in detail the combinatorial algorithm which yields the
$CM_{\mathbf{h}}$-partition of $\mathcal{P}(l,n)$. This is based on
\cite[$\S$ 7-8]{Go2} and is stated explicitly in \cite{GoMar}. Let $\Z^l_0 = \{ (\theta_0,
\dots,\theta_{l-1})\in \Z^l: \theta_0 + \dots + \theta_{l-1} = 0\}$,
an $S_l$ stable sublattice of $\Z^l$. There is an $S_l$-equivariant
isomorphism of lattices
\begin{align*}\phi : \Z_0^l \to R;\ (r_0, \dots ,r_{l-1}) \mapsto
(r_{l-1}-r_0)e_0 + (r_0-r_1)e_1 + \dots + (r_{l-2} -
r_{l-1})e_{l-1}.\end{align*} Let $S_l$ act on $\mathcal{P}(l,n)$ by $s_i(\lambda^{(0)}, \dots,
\lambda^{(l-1)})=(\lambda^{s_i(0)}, \dots, \lambda^{s_i(l-1)})$ for
all $i$.

\begin{thm}\cite[Theorem 2.5]{GoMar} Let $\mathbf{h} = (-1,H_0, H_1, \dots ,H_{l-1})\in \Q^{l+1}$ be a
parameter for $H_{\mathbf{h}}$ as in \ref{parameters}. Let $\theta =
(1+H_0, H_1, \dots, H_{l-1}) \in \Q_1^l.$ Let $w_{\theta} \in
\tilde{S}_l$ be such that $w_{\theta} \cdot \theta \in \mathcal{A}$
and suppose that $\theta$ has type $J$. Write $w_{\theta} =
\phi(\mathbf{r}) w$ with $w \in S_l$ and $\mathbf{r} \in \Z_0^l$.
Let $\blambda, \bmu \in \mathcal{P}(l,n)$. Then
\[\blambda \sim_{CM_{\mathbf{h}}} \bmu\ \mathrm{if\ and\ only\ if}\ \tau_{\mathbf{r}}(w(\blambda))_J =
\tau_{\mathbf{r}}(w(\bmu))_J.\]
\end{thm}

We will refer to $\theta$ as the stability parameter associated to
$\mathbf{h}$. Following this theorem we define the
\textit{$J$-heart} of a multipartition $\blambda$ to be the
partition $\tau_{\mathbf{r}}(w(\blambda))_J$. Note that this
definition depends on the parameter $\mathbf{h}$.

\subsection{}\label{Jthmsetup} Retain the hypotheses of Theorem
\ref{Iainthm}. Let $\bep:=\phi(\mathbf{r})w \cdot \theta =
(\varepsilon_0, \dots ,\varepsilon_{l-1})$, so that $0 \leq
\varepsilon_i \leq 1$ for all $i$ and $\varepsilon_0 + \dots +
\varepsilon_{l-1} = 1$. Recall the integer $d$ from
\ref{parameters}. A straightforward calculation shows that $d\bep
\in \Z^l$.

\begin{defn}
Let $\mathbf{c} = (c_0, \dots ,c_{l-1}) \in \Z^l$. Then for each $0
\leq i \leq l-1$ define the partial sum $m_i(\mathbf{c})= c_0 +
\dots + c_i$. %Define also $m_{-1}(\mathbf{c})=0$.
\end{defn}

We define an equivalence relation on the set $\{0, \dots ,l-1\}$
via: \begin{align}\label{equivreln} p\sim q\ \mathrm{if\ and\ only\
if}\ m_p(d\bep) = m_q(d\bep).
\end{align} For any $0\leq t \leq d$, let $I_t$ denote the equivalence class of
$p\in \{0,\dots ,l-1\}$ such that $m_p(d\bep) = t$. Given $0 \leq
a\leq b \leq l-1$, denote the corresponding interval by $[a,b]$. We
define $[a,b]=\emptyset$ if $a>b$.

\begin{lem}
Let $0 \leq t \leq d$. Then the set $I_t$ is empty or equal to
$[a,b]$ for some $0\leq a \leq b \leq l-1$. Furthermore, if $I_t \neq \emptyset$
then \[J \cap I_t = \begin{cases} [a+1,b]\ &\mathrm{if}\ t\neq0\\ [a,b]\ &\mathrm{if}\ t=0.\end{cases} \]
In particular, if $0\neq p \in J\cap I_t$ for some $t$ then $p-1 \in I_t$.
\end{lem}
\begin{proof} This follows immediately from the inequalities \begin{align*} 0\leq m_0(d\bep) \leq m_1(d\bep) \leq
\dots \leq m_{l-1}(d\bep) = d\end{align*} and the fact that, for all
$0\leq p \leq l-2$, $m_p(d\bep) = m_{p+1}(d\bep)$ if and only if
$\varepsilon_{p+1} = 0$.
\end{proof}

\subsection{}\label{Jheartdescrip} In order to calculate $\tau_{\mathbf{r}}(w(\blambda))_J$
one first considers, for each $0\leq p \leq l-1$, the
$\beta$-numbers
\[\beta^{r_p}(\lambda^{(w^{-1}(p))}) = \{ \lambda^{(w^{-1}(p))}_1 + r_{p},
\lambda^{(w^{-1}(p))}_2 + r_{p} - 1, \dots \}.\] Let $C^p =
C^p(\blambda) = \beta^{r_p}(\lambda^{(w^{-1}(p))})$ and let
$\mathbf{C} = (C^0, \dots ,C^{l-1})$. Let $0\leq t \leq d$ and let
$I_t$ be the corresponding equivalence class defined in
\ref{Jthmsetup}. For $1 \leq t \leq d-1$, let $C_{[t]}$ be the
multiset $\bigcup_{p\in I_t} C^p$. We define also the multiset
$C_{[0]} = C_{[d]} := \bigcup_{p\in I_0} S^{-1}C^p \cup
\bigcup_{p\in I_d} C^p$. We define $\pi(C_{[t]}):= \sum_{p\in I_t} \pi(C^p)$ for $1\leq t \leq d-1$ and
$\pi(C_{[0]}) = \pi(C_{[d]}) := \sum_{p\in I_0} x^{-1}\pi(C^p) + \sum_{p\in I_d} \pi(C^p)$.

We define $\tilde{\mathbf{C}}= \tilde{\mathbf{C}}(\blambda) = (\tilde{C}^0, \dots
,\tilde{C}^{l-1})$ to be the unique tuple of sets such that
\begin{itemize}\item[(i)] $\tilde{C}^p$ is a strictly decreasing set of
integers for all $p$,\item[(ii)] $\tilde{C}_{[t]} = C_{[t]}$ for all $0\leq t \leq d$,
\item[(iii)] if $0\neq p \in J$ then $C^p \subseteq C^{p-1}$ and
\item[(iv)] if $0\in I_0$ (so that $0\in J$) then $S^{-1}C^0 \subseteq C^{l-1}$. \end{itemize}

Such a tuple, $\tilde{\mathbf{C}}$, exists and is unique. Indeed, if
$\mathbf{C}$ does not satisfy (iii) then there is some $p$ and an
element $x\in C^p$ with $x\notin C^{p-1}$. We form a new tuple of
decreasing sequences by removing $x$ from $C^p$ and adding it to
$C^{p-1}$ - this new tuple still satisfies (ii) by Lemma
\ref{Jthmsetup}. We now repeat this process (and the analogous one
for (iv)) until we have a tuple with the desired properties. This
shows that there exists a $\tilde{\mathbf{C}}$ satisfying (i)-(iv).
Suppose that there exists another tuple $\mathbf{D}$ satisfying
these conditions and that $\mathbf{D} \neq \tilde{\mathbf{C}}$. Thus
there is some $p$ and some $i$ such that $D^p_i \notin
\tilde{C}^p_i$. Suppose that $p\in I_t$ for some $1\leq t \leq d-1$,
and let us set $x = D^p_i$. By Lemma \ref{Jthmsetup}, $I_t$ is of
the form $[a,b]$ for some $a$ and $b$. By (iii), $x\in D^p \subseteq
D^{p-1} \subseteq \dots \subseteq D^a$. By properties (ii) and
(iii), there exists a $k$ such that $x\in \tilde{C}^k \subseteq
\tilde{C}^{k-1} \subseteq \dots \subseteq \tilde{C}^a$ and $x\notin
\tilde{C}^{k+1}, \tilde{C}^{k+2}, \dots , \tilde{C}^b$. By our
assumption on $x$, $k < p$. This is a contradiction, since both
$\mathbf{D}$ and $\tilde{\mathbf{C}}$ satisfy (ii). A similar
argument when $p\in I_0 \cup I_d$ (using both conditions (iii) and
(iv)) also yields a contradiction. Therefore $\tilde{\mathbf{C}}$ is
uniquely defined.

\begin{thm}\label{Jheartthm} The $J$-heart of
$\tau_{\mathbf{r}}(w(\blambda))$ is the partition with
$\beta$-number equal to $\chi(\tilde{\mathbf{C}})$.
\end{thm}
\begin{proof} Let $\mu$ denote the partition with $\beta$-number
$\chi(\tilde{\mathbf{C}})$. We first show that $\mu_J = \mu$. Let $y
\in Y(\mu)$. Suppose that $y$ is $J$-removable. By \ref{Jhearts},
$y=(k,\mu_k)$ for some $k$. The element $\beta(\mu)_k \in
\chi(\tilde{\mathbf{C}})$ is equal to $l{\tilde{C}}^p_i+p-l+1$ for
some $0\leq p\leq l-1$ and $i\geq 1$. By \ref{Jhearts}, the
$J$-removability of $y$ is equivalent to the fact that there is a
$j\in J$ such that $l{\tilde{C}}^p_i+p-l+1 - 1 \equiv p \equiv j\
\mathrm{mod}\ l$ and $\beta(\mu)_{k+1} < l{\tilde{C}}^p_i+p-l$.

Let us first suppose that $0\neq p \in J$. Consider the equality
$\beta(\mu)_{k+1} = l{\tilde{C}}^p_i+p-l$. Since $\beta(\mu)_{k+1}
\equiv p$ mod $l$, this is equivalent to: $\beta(\mu)_{k+1}=
l{\tilde{C}}^{p-1}_{i'}+(p-1)-l+1$ for some $i'$ and
${\tilde{C}}^{p-1}_{i'} = {\tilde{C}}^p_i$. Thus the
$J$-removability of $y$ is equivalent to $\tilde{C}^p_i\notin
\tilde{C}^{p-1}$. By the construction of $\tilde{\mathbf{C}}$,
$\tilde{C}^p \subseteq \tilde{C}^{p-1}$ for all $0\neq p \in J$,
which is a contradiction.

Suppose now that $p=0\in J$ and so $I_0 \neq \emptyset$. Then an
analogous argument to that given in the previous paragraph shows
that the $J$-removability of $y$ is equivalent to ${\tilde{C}}^0_i-1 \notin
{\tilde{C}}^{l-1}$. Again we obtain a contradiction to our construction of
$\tilde{\mathbf{C}}$.

Let $\nu = \tau_{\mathbf{r}}(w(\blambda))$. We now show that we can
obtain $\mu$ from $\nu$ by removing $J$-removable nodes. If $\nu$
has no $J$-removable nodes then, as follows from the arguments in
the previous two paragraphs, this implies $C^p \subseteq C^{p-1}$
for all $0\neq p \in J$ and $S^{-1}C^0 \subseteq C^{l-1}$ if $0\in
J$. Thus $\mathbf{C}$ satisfies conditions (iii) and (iv) above.
Therefore $\mathbf{C}= \tilde{\mathbf{C}}$ and $\nu_J = \nu = \mu$.
If $\nu$ has a $J$-removable node then this means there is some
$0\neq p\in J$ (respectively $0=p\in J$) and some $i$ such that
$C^p_i \notin C^{p-1}$ (respectively $C^0_i-1 \notin C^{l-1}$). By
\ref{Jhearts}, the partition obtained by removing this node has
$\beta$-number $\chi(\mathbf{C}')$ where $\mathbf{C}'=({C'}^0, \dots
,{C'}^{l-1})$ with ${C'}^p = C^p\setminus\{C^p_i\}$, ${C'}^{p-1}=
C^{p-1} \cup \{C^p_i\}$ (respectively ${C'}^{l-1}= C^{l-1} \cup
\{C^0_i-1\}$), and ${C'}^q = C^q$ otherwise. Now $\chi(\mathbf{C}')
= \beta(\nu_J)$ if and only if $\mathbf{C}'$ satisfies the defining
properties of $\tilde{\mathbf{C}}$. If this is not the case then we
can repeat this process: after removing $k$ $J$-removable nodes we
have a partition with $\beta$-number $\chi(\mathbf{C}^{(k)})$ for an
$l$-tuple $\mathbf{C}^{(k)}$ satisfying (i) and (ii) above. This
partition has no $J$-removable boxes precisely when
$\mathbf{C}^{(k)}$ also satisfies conditions (iii) and (iv), that
is, $\mathbf{C}^{(k)}=\tilde{\mathbf{C}}$. Therefore $\beta (\nu_J)
= \chi(\tilde{\mathbf{C}}) = \beta(\mu)$ and so $\nu_J = \mu$.
\end{proof}

\begin{cor}
Two multipartitons $\blambda, \bmu \in \mathcal{P}(l,n)$ have the same
$J$-heart if and only if $C(\blambda)_{[t]} = C(\bmu)_{[t]}$ for all $0\leq t \leq d$.
\end{cor}
\begin{proof} If $\blambda$ and $\bmu$ have the same
$J$-heart then $\chi(\tilde{\mathbf{C}}(\blambda)) =
\chi(\tilde{\mathbf{C}}(\bmu))$. Thus $\tilde{\mathbf{C}}(\blambda)$
and $\tilde{\mathbf{C}}(\bmu)$ and in particular, $C(\blambda)_{[t]}
= C(\bmu)_{[t]}$ for all $0\leq t \leq d$. On the other hand, if
$C(\blambda)_{[t]} = C(\bmu)_{[t]}$ for all $0\leq t \leq d$, then
by the uniqueness of $\tilde{\mathbf{C}}(\blambda)$ and
$\tilde{\mathbf{C}}(\bmu)$ we have $\tilde{\mathbf{C}}(\blambda) =
\tilde{\mathbf{C}}(\bmu)$. Therefore $\blambda, \bmu$ have the same
$J$-heart by the theorem.
\end{proof}

\subsection{}\label{Heckeparameters}
Recall the action of $\tilde{S}_l$ on $\Z^l$ defined in
(\ref{affinesymaction}).

\begin{lem} Let $d$ be as in \ref{parameters} and let $\bep$ be as in \ref{Jthmsetup}.
Let $\theta$ be the stability parameter associated to $\mathbf{h}$, and let $w_{\theta} =
\phi(\mathbf{r}) w$ with $w \in S_l$ and $\mathbf{r} \in \Z_0^l$. Then
\begin{itemize}
\item[(1)] for all $0 \leq p \leq l-1$, \begin{align*}m_p(w \cdot
\theta) = m_{w^{-1}(p)}(\theta);\end{align*}
\item[(2)] let $0\leq t\leq d$, then for all $p\in I_t$, $m_{w^{-1}(p)}(d\theta) = dr_p
- dr_{l-1} + t$.
\end{itemize}
\end{lem}
\begin{proof}
Part (1) follows from the easy fact that $m_p(s_k \cdot
\theta) = m_{s_k(p)}(\theta)$ for all $1\leq k\leq l-1$. Part (2) follows from
\[m_{w^{-1}(p)}(d\theta) = m_{w^{-1}(p)}(w^{-1}(\phi(-d\mathbf{r})+ d\bep))= dm_{p}(\phi(-\mathbf{r}))+ m_{p}(d \bep).
\] and the equalities
$dm_{p}(\phi(-\mathbf{r}))=dr_p - dr_{l-1}$ and $m_{p}(d \bep) = t$.
\end{proof}

\subsection{} We can now prove our main theorem.

\begin{thm}\label{mainthm}
Let $\blambda, \bmu \in \mathcal{P}(l,n)$. Then $\blambda$ and
$\bmu$ have the same $J$-heart if and only if
$\mathrm{Res}^{\overline{\mathbf{m}}}_{\blambda}(x^d) =
\mathrm{Res}^{\overline{\mathbf{m}}}_{\bmu}(x^d)$.
\end{thm}
\begin{proof}
Let $K$ be a positive integer such that $C(\blambda)^p$ stabilises
with respect to $r_p$ after $K$ steps, for all $p$ and all $\blambda
\in \mathcal{P}(l,n)$ (we could, for instance, take any $K>n$). Let
$\blambda, \bmu \in \mathcal{P}(l,n)$. By Corollary
\ref{Jheartdescrip}, $\blambda$ and $\bmu$ have the same $J$-heart
if and only if $C(\blambda)_{[t]} = C(\bmu)_{[t]}$, or equivalently,
$\pi(C(\blambda)_{[t]}) = \pi(C(\bmu)_{[t]})$ for all $0\leq t \leq
d$.

Suppose that $1\leq t \leq d-1$. By our choice of $K$,
$\pi(C(\blambda)_{[t]}) = \pi(C(\bmu)_{[t]})$ if and only if
$\sum_{p\in I_t} \pi(C(\blambda)^p)_{\leq K} = \sum_{p\in I_t}
\pi(C(\bmu)^p)_{\leq K}$, and this latter equality is equivalent to
\begin{align}\label{equation1} \sum_{p\in I_t} \pi(C(\blambda)^p)_{\leq K}(x^d) = \sum_{p\in I_t}
\pi(C(\bmu)^p)_{\leq K}(x^d).\end{align} By (\ref{res-cont}),
\[\sum_{p\in I_t} \pi(C(\blambda)^p)_{\leq K}(x^d) =
(x^d-1)\sum_{p\in I_t} \Res^{dr_p}_{\lambda^{(w^{-1}(p))}}(x^d) +
\sum_{p\in I_t} \frac{x^{dr_p+d-dK}-x^{dr_p+d}}{1-x^d}.\] Thus
(\ref{equation1}) is equivalent to \begin{align}\label{equation2}
\sum_{p\in I_t} \Res^{dr_p}_{\lambda^{(w^{-1}(p))}}(x^d) =
\sum_{p\in I_t} \Res^{dr_p}_{\mu^{(w^{-1}(p))}}(x^d). \end{align} By
Lemma \ref{Heckeparameters}(2), (\ref{equation2}) is equivalent to
\begin{align*} x^{dr_{l-1}-t}\sum_{p\in I_t}
\Res^{m_{w^{-1}(p)}(d\theta)}_{\lambda^{(w^{-1}(p))}}(x^d) =
x^{dr_{l-1}-t}\sum_{p\in I_t}
\Res^{m_{w^{-1}(p)}(d\theta)}_{\mu^{(w^{-1}(p))}}(x^d),\end{align*}
which is equivalent to \begin{align}\label{equation3} \sum_{p\in
I_t} \Res^{m_{w^{-1}(p)}(d\theta)}_{\lambda^{(w^{-1}(p))}}(x^d) =
\sum_{p\in I_t}
\Res^{m_{w^{-1}(p)}(d\theta)}_{\mu^{(w^{-1}(p))}}(x^d).\end{align}

An analogous argument for $t=0$ yields $C(\blambda)_{[0]} =
C(\bmu)_{[0]}$ if and only if \begin{align}\label{equation4}
\sum_{p\in I_0 \cup I_d}
\Res^{m_{w^{-1}(p)}(d\theta)}_{\lambda^{(w^{-1}(p))}}(x^d) =
\sum_{p\in I_0 \cup I_d}
\Res^{m_{w^{-1}(p)}(d\theta)}_{\mu^{(w^{-1}(p))}}(x^d).
\end{align}

Let $m(d\theta) = (m_0(d\theta), \dots ,m_{l-1}(d\theta))$. We claim
that equalities (\ref{equation3}) for all $1\leq t \leq d-1$ and
(\ref{equation4}) together are equivalent to
$\mathrm{Res}^{m(d\theta)}_{\blambda}(x^d) =
\mathrm{Res}^{m(d\theta)}_{\bmu}(x^d)$. Consider the decomposition
of vector spaces $\Z[x, x^{-1}] = \bigoplus_{0\leq k \leq d-1} x^k
\Z[x^d, x^{-d}]$, and let $P_k$ denote the projection onto the $k$th
component. By Lemma \ref{Heckeparameters}(2),
\[P_k(\mathrm{Res}^{m(d\theta)}_{\blambda}(x^d)) = \begin{cases}\sum_{p\in I_k}
\Res^{m_{w^{-1}(p)}(d\theta)}_{\lambda^{(w^{-1}(p))}}(x^d)\
&\mathrm{if}\ 1\leq k\leq d-1\\ \sum_{p\in I_0 \cup I_d}
\Res^{m_{w^{-1}(p)}(d\theta)}_{\lambda^{(w^{-1}(p))}}(x^d)\
&\mathrm{if}\ k=0.
\end{cases} \] Since $\mathrm{Res}^{m(d\theta)}_{\blambda}(x^d) =
\mathrm{Res}^{m(d\theta)}_{\bmu}(x^d)$ if and only if
$P_k(\mathrm{Res}^{m(d\theta)}_{\blambda}(x^d)) =
P_k(\mathrm{Res}^{m(d\theta)}_{\bmu}(x^d))$ for all $k$, this proves
the claim.

To complete the proof of the theorem, we note that for all $\blambda
\in \mathcal{P}(l,n)$, $\mathrm{Res}^{m(d\theta)}_{\blambda}(x^d)=
x^{d+dH_0}\mathrm{Res}^{\overline{\mathbf{m}}}_{\blambda}(x^d).$
\end{proof}

We now verify Conjecture \ref{families}(i) for the groups
$G(l,1,n)$.

\begin{cor}
Let $\blambda, \bmu \in \mathcal{P}(l,n)$. Then $\blambda
\sim_{R_{\mathbf{m}}} \bmu$ implies $\blambda
\sim_{CM_{d\mathbf{h}}} \bmu$.
\end{cor}
\begin{proof}
By Theorem \ref{multiparts}, if $\blambda \sim_{R_{\mathbf{m}}}
\bmu$ then $\mathrm{Res}^{\overline{\mathbf{m}}}_{\blambda}(x^d) =
\mathrm{Res}^{\overline{\mathbf{m}}}_{\bmu}(x^d)$. The above theorem
then implies that $\blambda$ and $\bmu$ have the same $J$-heart and
therefore, by Theorem \ref{Iainthm}, $\blambda
\sim_{CM_{\mathbf{h}}} \bmu$. By Lemma \ref{CMpart}, $\blambda
\sim_{CM_{\mathbf{h}}} \bmu$ if and only if $\blambda
\sim_{CM_{d\mathbf{h}}} \bmu$.
\end{proof}

We conclude with a mention of another interpretation of the $CM_{\bf
h}$-partition for $G(l,1,n)$. Let us suppose that the parameter $\bf
h$ has integer entries. For the corresponding tuple, ${\mathbf{m}}$,
one can define a representation, $F(\Lambda_{\mathbf{m}})$, of the
quantum algebra $U_v(\mathfrak{sl}_{\infty})$, see \cite[Section
6]{LM} for further details. This module is called the higher level
Fock space, and it has a basis of weight vectors, $s_{\blambda}$,
labelled naturally by multipartitions of $n$. As a consequence of
Theorem \ref{mainthm} and \cite[Section 6.2]{LM} we have
\begin{cor}\label{Fockweights} Let $\blambda, \bmu \in \mathcal{P}(l,n)$. Then
$\blambda \sim_{CM_{\mathbf{h}}} \bmu$ if and only if $s_{\blambda}$
and $s_{\bmu}$ have the same weight.
\end{cor}

\bibliographystyle{alphanum}
\bibliography{families}

\end{document}